\documentclass[10pt]{article}
\usepackage{amsmath}
\usepackage{fullpage}
\usepackage{amsfonts}
\usepackage{amssymb}
\usepackage{amsthm}

\numberwithin{equation}{section}
\theoremstyle{plain}
\newtheorem{theorem}{Theorem}[section]
\newtheorem{lemma}{Lemma}[section]
\newtheorem{corollary}{Corollary}[section]

\begin{document}

\title{Algorithmically random Fourier series and Brownian motion}

\author{Paul Potgieter}
\date{}
\maketitle \vspace{-1.2cm}
\begin{center}
 \emph{Department of Decision Sciences, University
of South Africa\\ P.O. Box 392, Pretoria 0003, South Africa}\\
\texttt{potgip@unisa.ac.za}
\end{center} 

\begin{abstract}
We consider some random series parametrised by complex binary strings. The simplest case is that of Rademacher series, independent of a time parameter.  This is then extended to the case of Fourier series on the circle with Rademacher coefficients. Finally, a specific Fourier series which has coefficients determined by a computable function is shown to converge to an algorithmically random Brownian motion.

\vspace{0.5cm} \noindent {\bf Mathematics Subject Classification
(2010)} Primary 68Q30, 42A20; secondary 42A38, 60G15 

\end{abstract}

\section{Introduction} 
\label{intro}
There has been much success in recent years in interpreting probabilistic phenomena in the language of algorithmic randomness. Of particular interest to our current investigations are the effective versions of the Borel-Cantelli lemmas by Davie~\cite{Davie} and the construction of algorithmically random Brownian motion by Fouch\'{e}~\cite{Fouche}. Using these, it has been shown that many effectivised properties of Brownian motion hold for a version of Brownian motion parametrised by Kolmogorov-Chaitin complex strings. Important examples of such properties are the law of the iterated logarithm~\cite{Hanssen} and results on the Fourier spectra of measures on the images of fractals~\cite{FMD}.

The aim of the first part of the paper (sections 2 and 3) is to formulate convergence and divergence results for algorithmically random Rademacher series which are determined by KC-strings. This allows us to state with certainty when certain series will converge. The proofs rely on certain Kolmogorov and Paley-Zygmund inequalities in Hilbert space. Although only the real case is considered, the proofs should be transferable to other spaces where a suitable notion of computability is defined. We then turn to establishing similar results for trigonometric series with coefficients again determined by Rademacher series. In particular, we will show that under suitable conditions, the algorithmically random Rademacher series forms the Fourier coefficients of a Borel measure on $[0,1]$. As possible future work, it would be interesting to study the class of such measures in terms of computability. 

The second part of the paper (sections 4 and 5) is concerned with the convergence of computably determined series to continuous functions. This is done firstly for algorithmically random Rademacher series. We then consider a particular series, the Fourier-Wiener series, used in a construction of Brownian motion, and show that this not only converges to a continuous function, but to a complex oscillation in the sense of Asarin and Pokrovskii. This provides an alternative to Fouch\'{e}'s construction of algorithmically random Brownian motion~\cite{Fouche}. 

The probabilistic aspects of this paper are heavily influenced by the exposition of Kahane~\cite{Kahane}. For computability aspects, we shall throughout notation similar to that used in~\cite{Fouche}, which we summarise here.

We denote the set of non-negative integers by $\omega$ and write $\mathcal{N}$ for the space $\{-1,1\}^{\omega}$. The subset of finite strings is denoted by $\{-1,1\}^{*}$.  For $\alpha = \alpha_0 \alpha_1 \dots$ in $\mathcal{N}$, we write $\overline{\alpha}(n) = \prod_{j<n} \alpha_j$ for the intitial segment of $\alpha$ of length $n$. For $a\in \{-1,1\}^{*}$, we write $|a|$ for the length of a and $[a]$ for the interval generated by $a$; that is, the set $\{\alpha \in \mathcal{N}: \overline{\alpha}(|a|) = a\}$.  
For $\alpha \in \mathcal{N}$, we denote the Kolmogorov complexity of $\alpha$ by $K(\alpha)$, and say that $\alpha$ is Kolmogorov-Chaitin complex if $\exists k \forall n K(\overline{\alpha}(n))\geq n-k$. The set of Kolmogorov-Chaitin complex strings shall be denoted by $KC$.

A sequence $(a_n)$ of real numbers is said to converge to $0$ effectively if there is some computable function $f:\omega \to \omega$ such that $|a_n |\leq (m+1)^{-1}$ whenever $n\geq f(m)$, for all $n,m \in \omega$. Letting $\lambda$ be the Lebesgue probability measure over $\mathcal{N}$, we say that a subset $A$ of $\mathcal{N}$ has constructive measure $0$ if there is a computable function $g:\omega^2 \to \{-1,1\}^{*}$ such that $A \subset \cap_n \cup_m [g(n,m)]$, where $\lambda (\cup_m [g(n,m)])$ converges effectively to $0$ as $n\to \infty$.  

We shall use the following theorems or their probabilistic analogues repeatedly. 

\begin{theorem}\label{thm1.1} \cite{MartinLof}
If $\alpha \in KC$, $\alpha$ is in the complement of every set of constructive measure $0$.
\end{theorem}

\begin{theorem}\label{thm1.2} \cite{Fouche3}
If $B$ is a $\Sigma^{0}_{1}$ set of Lebesgue measure $1$, then $B$ contains $KC$.
\end{theorem}

This theorem has the converse:

\begin{theorem}\label{thm1.3}
A $\Pi_{1}^{0}$ set of Lebesgue measure $0$ contains no element of $KC$.
\end{theorem}

The following is an effective version of the Borel-Cantelli lemmas. 

\begin{lemma}\label{lem1.1}~\cite{Davie} Let $\{A_k \}$ be a uniform sequence of $\Sigma^{0}_{1}$ sets.
\begin{enumerate}
\item[(i)]{If $\sum_k \lambda (A_k) < \infty$, then, for each $\alpha \in KC$, it is this case that $\alpha$ belongs to no $A_k$ for sufficiently large $k$.}
\item[(ii)]{If the events $\{ A_k \}$ are pairwise independent and $\sum_k \lambda (A_k) = \infty$, then $A_k$ must occur infinitely often for each $\alpha \in KC$.}
\end{enumerate}
\end{lemma}

Corresponding notions for Wiener measure will not be needed until section 5, and will only be introduced then.

\section{Algorithmically random Rademacher series}

The goal of this section is to establish convergence and divergence results for algorithmically random Rademacher series. Since the convergence of a Rademacher series is a tail phenomenon, a given series converges almost surely, or diverges almost surely. In this section we will show that if a Rademacher series of a computable sequence of numbers converges a.s., it will converge for each KC-string as long as  the $l^2$ norm of the sequence is a computable number. Conversely, we show that if a Rademacher series of a computable sequence of numbers diverges, it diverges for each KC-string. A consequence of this is that convergence/divergence for all KC-strings can be determined by the convergence/divergence of any one such string. 

A Rademacher sequence on a probability space $(\Omega, \Sigma ,\mathbb{P})$ is defined as a sequence  $(\varepsilon_n)$ of random variables on $X$ such that $\mathbb{P}(\varepsilon_n = 1)=\mathbb{P}(\varepsilon_n = -1)=1/2$, $n=1,2,\dots$. In this paper, we examine Rademacher sequences parametrised by $\mathcal{N}$, which is to say that for some $x\in \mathcal{N}$, $\varepsilon_n (x) =-1$ or $1$ according to whether the $n$th bit of $x$ is $-1$ or $1$. A Rademacher series in a Hilbert (or Banach) space is a series of the form $\sum_{n=1}^{\infty} \varepsilon_n u_n $, where $(\varepsilon_n)$ is a Rademacher sequence and $(u_n)$ is a sequence of vectors in a Hilbert (or Banach) space. We shall only concern ourselves with real sequences $(u_n)$. A sequence $(\varepsilon_n)$ is called an \emph{algorithmically random} Rademacher sequence if there is some $\alpha \in KC$ such that $\varepsilon_k$ is the $k$th bit of $\alpha$. When it is necessary to specify $\alpha$, the sequence is denoted by $(\varepsilon_n (\alpha))$. 

A sequence $(x_n)$ of real numbers is said to be a \emph{computable sequence} if there is an effective way of approximating $x_n$ arbitrarily closely, given $n$. We shall say that a computable sequence of real numbers $(x_n)$ is \emph{square computable} if there is an effective procedure for determining the sum $\sum_{k=1}^{\infty} x_{k}^{2}$ to arbitrary accuracy from $(x_n)$. In other words, given a positive integer $m$, there is an effective procedure for finding $n=n(m)$ such that 
\begin{equation}
\left| \sum_{k=1}^{\infty} x_{k}^{2} - \sum_{k=1}^{n} x_{k}^{2} \right| < 1/m.
\end{equation}

Throughout, $\mathbb{E}(X)$ shall denote the expected value of a random variable $X$, and $V(X)= \mathbb{E}((X-\mathbb{E}(X))^2)$ its variance. All events considered are members of the usual Borel $\sigma$-algebra over $\mathcal{N}$, and the probability measure is denoted by $\mathbb{P}$ rather than $\lambda$, in order to maintain a more uniform notation.

The first relevant theorem is the following. For generality, we phrase the theorem with the sequence $( u_n )$ considered as a sequence in a Hilbert space $H$, with $\| \cdot \|$ denoting the norm associated with the inner product on $H$. 

\begin{theorem}\label{thm2.1}(\cite{Kahane} p30) Suppose that $X_n \in L_{H}^{2}(\Omega)$ are independent random variables and $\mathbb{E}(X_n )= 0$ for each $n$, and $\sum_{1}^{\infty} V(X_n) < \infty$. Then the series $\sum_{1}^{\infty}X_n$ converges in $H$ a.s. 
\end{theorem}

In particular, if our variables are $X_n = \varepsilon_n u_n$ for a Rademacher sequence $(\varepsilon_n)$ and if $\sum_{1}^{\infty} \|u_n \|^2 <\infty $, then $\sum_{1}^{\infty} \varepsilon_n u_n$ converges a.s. In the sequel, we shall omit the use of subscripts and superscripts in the summation unless confusion is likely. 

In order to prove a version of the theorem for complex strings, we shall need an inequality of Kolmogorov, which is also used in the proof of the Theorem \ref{thm2.1}.

\begin{theorem}\label{thm2.2}(\cite{Kahane} p29) Supposing that  $X_n \in L_{H}^2 (\Omega)$ are independent random variables and $\mathbb{E}(X_n)=0$ for each $n$, we have for any $N$
\begin{equation}\label{eq2.1}
\mathbb{P}\left( \sup_{n=1,2,\dots ,N} \| X_1 +X_2 +\cdots +X_n \| >r \right) < \frac{1}{r^2}\left( V(X_1) +V(X_2) + \cdots +V(X_N)\right).
\end{equation}
\end{theorem}

We can now prove a version of Theorem \ref{thm2.1} for complex strings.

\begin{theorem}\label{thm2.3}
Suppose that $(u_n)$ is a square computable sequence in $\mathbb{R}$, and $(\varepsilon_n )$ is an algorithmically random Rademacher sequence. Then $\sum \varepsilon_n u_n$ converges in $\mathbb{R}$.
\end{theorem}
\textbf{Proof.} Consider the event
\begin{equation}
B_{n}^{(k)} = \left( \sup_{j\geq 1} \left| \sum_{l=n}^{n+j} \varepsilon_l u_l \right| >\frac{1}{k} \right)
\end{equation}
for some positive integer $k$. By equation (\ref{eq2.1}), 
\begin{equation}
\mathbb{P}(B_{n}^{(k)}) < k^2 \sum_{n}^{\infty} u_{n}^{2}.
\end{equation}
Since we require that $( u_n )$ is square computable, we can effectively determine, for each $m$, an $n=n(m)$ such that $k^2 \sum_{n}^{\infty} u_{n}^{2}$ is smaller than $1/m^2$. 

Note that we can describe $B_{n}^{(k)}$ as
\begin{equation}
\bigcup_{j=1}^{\infty}\left( \left| \sum_{l=n}^{n+j} \varepsilon_l u_l \right|  >\frac{1}{k}\right).
\end{equation}
Thus $B^k_n$ is a countable union of finitely determined events, and is therefore a $\sum_{1}^{0}$ set. Similarly, since $n(m)$ can be effectively determined, the sequence $\{ B^{(k)}_{n(m)}\}_{m=1}^{\infty}$ forms a uniform sequence of $\sum_{1}^{0}$ sets. Moreover, the sum of the measures of the sets converges. By the effective version of the Borel-Cantelli lemma, we have that for a complex string $\alpha$ that $\alpha \notin B^{(k)}_{n(m)}$ for all large values of $m$. Therefore, given any positive $k\in \mathbb{Z}$, we have that, for large enough values of $n$,
\begin{equation}
\sup_{j}\left| \sum_{l=n}^{n+j} \varepsilon_l u_l \right|  \leq \frac{1}{k}. 
\end{equation}
We can conclude that the series converges for each complex string.\qed

The converse of Theorem \ref{thm2.1} for a general Hilbert space is the following.

\begin{theorem}\label{thm2.4} (\cite{Kahane} p31) Suppose that $\| X_n \| \in L^4 (\Omega)$ are independent random variables, $\mathbb{E}(X_n)=0$ and $\mathbb{E}(\| X_n \|^4 )\leq C V^2 (X_n)$ for all $n$, for some absolute constant $C>0$. If the series $\sum X_n$ is a.s. bounded, then $\sum_{1}^{\infty} V(X_n) < \infty$.
\end{theorem}
This has the trivial corollary:
\begin{corollary} For any sequence of real numbers $( u_n)$ (computable or not), if $\sum \varepsilon_n u_n$ is bounded for each KC string, then $\sum u_{n}^{2} < \infty$.
\end{corollary}

When the sum $\sum_{1}^{\infty} u_{n}^{2}$ diverges, the Rademacher series $\sum \varepsilon_n u_n$ will almost surely diverge (\cite{Kahane}, p31). For complex strings, we obtain the following result.

\begin{theorem}\label{thm2.5} If $( u_n )$ is a computable sequence such that $\sum u_{n}^{2} =\infty$, then $\sum \varepsilon_n u_n$ diverges for each algorithmically random Rademacher sequence $( \varepsilon_n )$.
\end{theorem}

We shall need the following Paley-Zygmund inequality for the proof.

\begin{lemma}\label{lem2.1}(\cite{Kahane} p31) If $( u_n)$ is a sequence of vectors in a Hilbert space and $( \varepsilon_n )$ is a Rademacher sequence, then for any $0<\lambda <1$ we have
\begin{equation}\label{eqlem2.1}
\mathbb{P}[\| \varepsilon_1 u_1 +\cdots +\varepsilon_m u_m \| > \lambda ( \|u_1 \|^2 +\cdots +\|u_m \|^2 )^{\frac{1}{2}}] >\frac{1}{3}(1-\lambda^2 )^2.
\end{equation}
\end{lemma}

\textbf{Proof of Theorem \ref{thm2.5}.} We can interpret Lemma \ref{lem2.1} as follows: For any positive integer $m$,
\begin{equation}
\mathbb{P}[ | \varepsilon_1 u_1 +\cdots +\varepsilon_m u_m |> \lambda ]  >\frac{1}{3}(1-\gamma^2)^2,
\end{equation}
where 
\begin{equation}
\gamma = \frac{\lambda}{(u_{1}^{2} +\cdots +u_{m}^{2})^{\frac{1}{2}}}.
\end{equation} 
In order for the inequality to apply, we need that $0<\gamma <1$. Since $\sum u_{n}^{2}$ diverges, we can take $m$ large enough in the above so that this is true for a given value of $\lambda$. Moreover, if $\lambda$ is computable, the smallest such $m$ can be effectively determined. For ease of computation, we set $\lambda = 1/2$. We can now effectively find $m_1$ so that the associated value of $\gamma$ is between $0$ and $1$ and
\begin{equation}
\eta_{m_1} = \frac{1}{3}(1-\gamma_{m_1}^{2})^2 = \frac{1}{3}\left( 1-\frac{\lambda^2}{u_{1}^{2} +\cdots +u_{m_1}^{2}} \right) >\frac{1}{6}.
\end{equation}
Since the inequality (\ref{eqlem2.1}) applies to all finite subsequences of $(u_n )$, given $m_k$ we can effectively find $m_{k+1} >m_k$ so that
\begin{equation}
\eta_{m_{k+1}} = \frac{1}{3}(1-\gamma_{m_{k+1}}^{2})^2 = \frac{1}{3}\left( 1-\frac{\lambda^2}{u_{m_k +1}^{2} +\cdots +u_{m_{k+1}}^{2}}\right) >\frac{1}{6}.
\end{equation}
Continuing in this way, we can effectively determine $m_{n+k}$, $k=1,2,\dots$. 
The events $(| \varepsilon_{m_k +1} u_{m_k +1} +\cdots +\varepsilon_{m_{k+1}} u_{m_{k+1}} | > \lambda )$, $k=1,2,\dots$, are independent finitely determined events, and moreover have probabilities that sum to $\infty$. This means, by the second effective Borel-Cantelli lemma, that they must occur infinitely often for complex strings. This implies that the sum for a complex string will always have sections of the tail which are larger than $1/2$ in absolute value, implying divergence. \qed
 
\begin{corollary}
If $(u_n )$ is a square computable sequence in $\mathbb{R}$ and $\sum \varepsilon_n u_n$ converges for any $\alpha \in KC$, then $\sum u_n^2$ converges. If $\sum \varepsilon_n u_n$ diverges for any $\alpha \in KC$, $\sum u_n ^2$ must diverge or converge to an incomputable number. 
\end{corollary}
\textbf{Proof.} Since the series $\sum u_{n}^{2}$ has only positive terms, it must converge or diverge to $\infty$. However, divergence implies that $\sum \varepsilon_n u_n$ also diverges. The second statement follows similarly. \qed

We end the section with a principle of contraction for KC-strings. The original can be found in \cite{Kahane}, p20.

\begin{theorem}\label{thm2.6}
If $( u_n)$ is a square computable sequence, then for a bounded computable sequence $( \lambda_n )$ with $\sup_n |\lambda_n |$ computable, the series $\sum \varepsilon_n \lambda_n u_n$ converges for each algorithmically random Rademacher sequence $( \varepsilon_n)$.
\end{theorem}
\textbf{Proof.} The proof follows that of Theorem \ref{thm2.3}, by replacing $\sum |u_n|^2$ by $(\sup |\lambda_n |) \sum u_{n}^{2}$. \qed 

\section{Algorithmically random Fourier series}

In this section, we consider algorithmically random versions of series of the form 
\begin{equation}\label{eq3.1}
\sum_{n=0}^{\infty} X_n \cos (2\pi nt +\Phi_n),
\end{equation}
for independent symmetric random variables $X_n e^{2\pi i\Phi_n}$, $n\geq 0$, and $t\in [0,1]$. In particular, we let $X_n = \varepsilon_n x_n$, where $x_n $ and $\phi_n = \Phi_n$ are fixed real numbers and $(\varepsilon_n)$ a Rademacher sequence parametrised by $\mathcal{N}$, for $n\geq 0$. Our purpose is to establish conditions on the various parameters in order to obtain convergence or divergence for series 
\begin{equation}\label{eq3.2}
\sum_{n=0}^{\infty} \varepsilon_n (\alpha) x_n \cos (2\pi nt+\phi_n)
\end{equation}
where $(\varepsilon_n (\alpha))$ is an algorithmically random Rademacher sequence.
 A discussion of the general case can be found in Chapter 5 of \cite{Kahane}.

A trigonometric series
\begin{equation}\label{eq3.3}
\sum_{n=0}^{\infty} x_n \cos (2\pi nt+\phi_n)
\end{equation}
with $x_n$, $\phi_n \in \mathbb{R}$, $n\geq 0$ is said to be a Fourier-Stieltjes series if there exists a measure $\mu$ on $[0,1]$ such that
\begin{equation}
x_0 \cos 2\pi \phi_0 = \int d\mu (t)
\end{equation}
and
\begin{equation}
x_n e^{\pm 2\pi i\phi_n} = 2 \int e^{\pm 2\pi int}d\mu (t), \quad n \geq 1.
\end{equation}
We say that ($\ref{eq3.3}) \in M$ if it is a Fourier-Stieltjes series. 

The Fej\'{e}r sums of (\ref{eq3.1}) are given by
\begin{equation}\label{eq3.4}
\sigma_N (t, \omega) = \sum_{n=0}^{N}\left( 1-\frac{n}{N} \right) X_n (\omega) \cos (2\pi nt+\Phi_n (\omega)), \quad N=1,2,\dots
\end{equation}
If the Fej\'{e}r sum of a series diverges, the series itself will diverge, although the converse is not necessarily true.
We shall denote the Fej\'{e}r sum by $\sigma_N (t)$, since the involvement of $\omega$ will usually be clear.  
Letting $\| \cdot \|_1$ denote the $L^1 ([0,1])$ norm, we have that (p136, \cite{Zygmund})
\begin{equation}\label{eq3.5}
(\ref{eq3.3}) \in M \iff \sup_N \| \sigma_N \|_1 < \infty.
\end{equation}
By Theorem 1 on p13 of \cite{Kahane}, we can conclude that
\begin{equation}\label{eq3.6}
(\ref{eq3.1}) \in M \textrm{ a.s. } \iff  \sup_N \| \sigma_N \|_1 < \infty \textrm{ a.s. }
\end{equation}

\begin{theorem}\label{thm3.1}
Let $( x_n )$, $(\phi _n)$ be computable sequences of real numbers so that $\sum x^{2}_{n} = \infty$ and let $( \varepsilon_n (\alpha ))$ be an algorithmically random Rademacher sequence. Then, 
\begin{equation}
\limsup_{N\to \infty} \int_{0}^{1}|\sigma_N (t)|dt = \infty.
 \end{equation}
\end{theorem}
Thus, given a computable sequence $(x_n)$ such that $\sum x_{n}^{2}$ diverges, the series (\ref{eq3.2}) will not be a Fourier-Stieltjes series if $\alpha \in KC$. This partially answers the question of why it is difficult to explicitly construct a series of the form $\sum \pm x_n \cos (2\pi nt+\varphi_n)$ which is not a Fourier-Stieltjes series when $\sum x_{n}^{2} =\infty$, even though it happens almost surely.

\bigskip

\textbf{Proof of Theorem \ref{thm3.1}.} Each of the $\sigma_N$ is a continuous function over a closed and bounded interval and we can approximate the integral by a Riemann sum. For a given $M$, we subdivide the interval $[0,1]$ into intervals of length $1/M$. Consider the event 
\begin{equation}
B_K = \exists N \exists M \left( \frac{1}{M}\sum_{k=0}^{M-1} |\sigma_{N}(k/M) |>K \right),
\end{equation} where $K$ is a positive integer. By Proposition 6 on p50 of \cite{Kahane}, each $B_K$ has probability $1$. What is more, each $B_K$ is a $\Sigma_{1}^{0}$ set, and therefore must contain all of $KC$.\qed

The following lemma tells us something of pointwise divergence. 

\begin{lemma}\label{lem3.2}
Let $( x_n )$, $(\phi _n )$ be computable sequences so that $\sum x^{2}_{n} = \infty$ and let $( \varepsilon_n (\alpha ))$ be an algorithmically random Rademacher sequence. Then the series (\ref{eq3.2}) diverges for almost every $t\in [0,1]$. 
\end{lemma}
\textbf{Proof.} The key to the proof is to adapt the second proof of Proposition 6, used in the previous theorem. There it is proved that 
\begin{equation}
\limsup_{N\to \infty}\int_{0}^{1} |\sigma_N (t)|dt = \infty,
\end{equation}
almost surely. Since the Paley-Zygmund inequality used in the proof is valid for all $t$ and $N$, and since by Proposition 4 on p49 of~\cite{Kahane} the series
\begin{equation}
\sum_{1}^{\infty}x_n^2 \cos^2 (2\pi nt+\phi_n)
\end{equation} diverges a.e., it is permissible to consider any closed interval of $[0,1]$ in the proof, rather than the interval itself. If $I=[a,b]$ is such an interval, the proof then implies that 
\begin{equation}
\limsup_{N\to \infty}\int_I |\sigma_N (t)|dt = \infty,
\end{equation}
almost surely. Letting $I=[k2^{-j}, (k+1)2^{-j}]$, $j\in \mathbb{Z}^{+}$, $k\in \{0,1,2, \dots ,2^j -1 \}$ be any dyadic subinterval of $[0,1]$, and considering the Riemann sum as in the previous proof, this means that the event
\begin{equation}
B_{K,j,k} = \exists N \exists M \left[ \frac{1}{M}\sum_{n=0}^{M-1} \left| \sigma_N \left( \left( k+\frac{n}{M}\right) 2^{-j} \right) \right|  >K \right]
\end{equation}
has measure $1$. Since it is $\Sigma_{1}^{0}$, it must hold for every $\alpha \in KC$. Hence, the integral of $|\sigma_N (t)|$ diverges over arbitrarily small intervals as $N\to \infty$, which implies that the Fej\'{e}r sum and hence the series diverges a.e.\qed

\medskip

We now turn to the subject of convergence. Convergence almost everywhere on $[0,1]$ when $\sum x_{n}^{2} < \infty$ is actually the least interesting case. For fixed $\alpha \in KC$, the coefficients form a square summable sequence and hence the Fourier series converges in norm to an element of $L^2 ([0,1])$. The 
Fej\'{e}r sums then have to converge to a function in $L^1([0,1])$, which means that $\sup_N \|\sigma_N \|_1$ must be bounded, implying that (\ref{eq3.2}) $\in M$. Imposing certain constructive constraints on the measures in question might lead to results of greater interest. However, we can say something about the convergence at computable numbers.

\begin{theorem}\label{thm3.2}
If $(x_n )$ is a square computable sequence and $(\phi_n)$ is a computable sequence, the series (\ref{eq3.2}) converges at each computable number for each  $\alpha \in KC$. 
\end{theorem}
\textbf{Proof.}  Noting that $( x_n \cos (2\pi nt+\phi_n))$ is a computable sequence such that the absolute value of the $n$th term is bounded by $|x_n|$, the proof of Theorem 2.3 can be applied directly.\qed

\section{Series as continuous functions}

In our construction of algorithmically random Brownian motion, it will be key to show that certain Fourier series converge to continuous functions. The notation and method of proof of the following theorem for algorithmically random Rademacher Fourier series will be important. 

For random variables $X_n$, each with finite variance and expectation $0$, we set
\begin{equation}\label{eq4.1}
s_j = \left( \sum_{2^j \leq n <2^{j+1}} \mathbb{E}(X^2_n) \right)^{\frac{1}{2}}, \quad j=0,1,2,\dots
\end{equation}
In the case to be considered below, we have that $X_n = \varepsilon_n x_n$, where the $x_n$ are now fixed real numbers, and  $(\varepsilon_n)$ is a Rademacher sequence parametrised by $\mathcal{N}$, and
\begin{equation}
s_j = \left( \sum_{2^j \leq n <2^{j+1}} x^2_n \right)^{\frac{1}{2}}, \quad j=0,1,2,\dots
\end{equation}
\begin{theorem}\label{thm4.1}
Whenever $X_n = \varepsilon_n x_n$, where $(\varepsilon_n)$ is an algorithmically random Rademacher sequence and $(x_n)$, $(\phi_n)$ are computable sequences, $(s_j)$ is a decreasing sequence and $\sum_{1}^{\infty} s_j$ is finite, the Fourier series (\ref{eq3.2}) converges to a continuous function.
\end{theorem}
\textbf{Proof.} As in~\cite{Kahane}, we set
\begin{equation}
P_k (t) = \sum_{N_k}^{N_{k+1}-1}\varepsilon_n x_n \cos (2\pi nt+\phi_n),
\end{equation}
where $N_k = 2^{2^k}$. 
Let $B_k$ be the event
\begin{equation}
\| P_k \|_{\infty} > 6 \left( \log N_{k+1} \sum_{N_k}^{N_{k+1}-1} x_n^2 \right)^{\frac{1}{2}}
\end{equation}
By theorem 2 on p69 of~\cite{Kahane}, presented in this paper as Theorem \ref{thm5.1}, 
\begin{equation}
\mathbb{P} (B_k) \leq \frac{8\pi}{N_{k+1}^2}.
\end{equation}
The theorem is applicable because $(X_n)$ is a subnormal sequence, that is, the variables are independent and satisfy $\mathbb{E} (e^{\lambda X_n}) \leq e^{\lambda^2 /2}$ for any $\lambda \in \mathbb{R}$.
Since each $P_k$ is a continuous computable function, we can state the lower bound on the uniform norm in terms of rational numbers. Setting 
\begin{equation}
F(k) = \left( \log N_{k+1} \sum_{N_k}^{N_{k+1}-1} x_n^2 \right)^{\frac{1}{2}},
\end{equation}
we can describe the event $B_k$ as
\begin{equation}
\exists M \exists n\leq M \left( \left| P_k \left( \frac{n}{M} \right) \right| > 6F(k) \right).
\end{equation}
The events $(B_k)$ therefore form a uniform sequence of $\Sigma_1^0$ sets, and the sum of all their probabilities is clearly finite. Therefore, for any algorithmically random Rademacher sequence $(\varepsilon_n)$, 
\begin{equation}
\| P_k \|_{\infty} \leq 6\left( 2^{k/2}\left( \sum_{2^k}^{2^{k+1}-1}s_j^2 \right)^{\frac{1}{2}}\right)
\end{equation}
for large $k$. If we now require that $\sum_{j=1}^{\infty}s_j <\infty$ and that the sequence $s_j$ is decreasing, we have that
\begin{equation}
\sum_{k=1}^{\infty} 2^{k/2}\left( \sum_{2^k}^{2^{k+1}-1}s_j^2 \right)^{\frac{1}{2}} <\infty.
\end{equation}
Therefore, $\sum_{k=1}^{\infty}\| P_k \|_{\infty}< \infty$, and the series $\sum_{k=1}^{\infty}P_k(t)$ must be uniformly convergent. \qed

It would now be natural to ask what would happen in a more general case where, instead of a Rademacher sequence, one has a sequence of independent, identically distributed random variables (for which a notion of computability is suitably defined) as coefficients. Although we cannot yet formulate a general result, the next section provides a specific instance of convergence.

\section{Algorithmically random Brownian motion}

We require some background in order to phrase analogues of Theorem \ref{thm1.1} and Lemma \ref{lem1.1} for Wiener measure. For a more complete introduction, see \cite{Fouche} and \cite{Fouche2}. Throughout, $\mathbb{W}$ shall denote Wiener measure on $\Sigma$, the Borel $\sigma$-algebra over $C[0,1]$ (topologised by the uniform norm).

Let $C_n$, $n\geq 1$ be the class of elements of $C[0,1]$ that vanish at $0$ and are linear with slope $\pm \sqrt{n}$ on the intervals $[(i-1)/n,i/n]$, $i=1,2,\dots n$. It is clear that each $x\in C_n$ can be encoded by a binary string $c(x)=a_1 \cdots a_n$ in $\{-1,1\}^{*}$ by setting $a_i =1$ if the slope is positive on $[(i-1)/n,i/n]$ and $a_i=-1$ otherwise. A sequence $(x_n)$ in $C[0,1]$ is \emph{complex} if for each $n$, $x_n \in C_n$ and there exists a constant $d>0$ (independent of $n$) such that $K(c(x_n))\geq n-d$ for all $n$. A function $x\in C[0,1]$ is a \emph{complex oscillation} if there is a complex sequence $(x_n)$ such that $\|x-x_n \|_{\infty}$ converges effectively to $0$ as $n\to \infty$. The collection of complex oscillations will be denoted by $\mathcal{C}$.

The next notion needed is that of an \emph{effective generating sequence}, formulated by Fouch\'{e} in~\cite{Fouche2}. Since we shall only use one example of such in this section, we do not present the general definition, only an instance thereof. Let $\mathcal{G}_0$ be a collection of sets in $\Sigma$ of the form 
\begin{equation}\label{eq5.10}
a_1 X(t_1) +\cdots +a_nX(t_n) \leq L \quad \textrm{or} \quad a_1 X(t_1) +\cdots +a_nX(t_n) < L
\end{equation}
where $X$ is a one-dimensional Brownian motion on $[0,1]$, the $a_j$, $0\leq t_j \leq 1$ are rational numbers and $L$ is a computable real. Let $(G_i: i<\omega)$ be an enumeration of $\mathcal{G}_0$ such that for any given $i$, the sign, denominator and numerator of each of the $a_j$, $t_j$ in the representation (\ref{eq5.10}) can be effectively computed, and the real $L$ can be approximated arbitrarily closely. 

For $G\in \mathcal{G}_0$, the open $\varepsilon$-neighbourhood $O_{\varepsilon}(G)$ of $G$ as given by (\ref{eq5.10}) is described by
\begin{equation}
a_1 X(t_1) +\cdots +a_nX(t_n) < L+\varepsilon \sum_j |a_j|.
\end{equation}
The neighbourhood $O_{\varepsilon}(G^c)$ of $G^c$ is described by
\begin{equation}
a_1 X(t_1) +\cdots +a_nX(t_n) > L-\varepsilon \sum_j |a_j|.
\end{equation}
It is shown in \cite{Fouche} that $\mathcal{G}_0 = (G_i :i<\omega)$ satisfies the following conditions, and hence forms an effective generating sequence:
\begin{enumerate}
\item[1.]{For $G\in \mathcal{G}_0$, we have for $F=O_{\varepsilon}(G)$, $F=O_{\varepsilon}(G^c)$, $F=G$ or $F=G^c$ that $\mathbb{W}(F) = \mathbb{W}(\overline{F})$.}
\item[2.]{There is an effective procedure that yields, for each sequence $0\leq i_i < \cdots < i_n <\omega$ and $k<\omega$, a binary rational $\beta_k$ such that
\begin{equation}
|\mathbb{W}(G_{i_1}\cap \cdots \cap G_{i_n}) - \beta_k | < 2^{-k}.
\end{equation}In other words, the Wiener measure of finite intersections of elements of $\mathcal{G}_0$ is computable.}
\item[3.]{For $n,i<\omega$,  a rational number $\varepsilon>0$ and $x\in C_n$, both the relations $x\in O_{\varepsilon}(G_n)$ and $x\in O_{\varepsilon}(G_{n}^{c})$ are computable in $x,\varepsilon ,i$ and $n$.}
\end{enumerate}

From an effective generating sequence $\mathcal{G}_0$ we can generate an algebra $\mathcal{G}$ in an effective way, i.e. there will be an enumeration $(T_i :i<\omega)$ of $\mathcal{G}$ such that, for given $i$, we can effectively describe $T_i$ as a finite union of finite intersections of elements of $\mathcal{G}_0$ or their complements. The sequence $(T_i :i<\omega)$ is called a \emph{computable enumeration} of $\mathcal{G}$, and we can refer to $\mathcal{G}$ as an \emph{effectively generated algebra}. 

We now present analogues in $\Sigma$ of the arithmetical sets we have used in the previous sections, and for which versions of Theorems \ref{thm1.1}, \ref{thm1.2} and Lemma \ref{lem1.1} will hold.
  
A sequence $(A_n)$ of sets in $\mathcal{G}$ is $\mathcal{G}$-computably enumerable if it is of the form $(A_n = S_{\phi(n)})$ for some total computable function $\phi:\omega \to \omega$ and $(S_i)$ a computable enumeration of $\mathcal{G}$. The union $\cup_n A_n$ is then a $\Sigma_{1}^{0}(\mathcal{G})$ set. The complement of such is a $\Pi_{1}^{0} (\mathcal{G})$ set. A sequence $(B_n)$ in $\mathcal{G}$ is a \emph{uniform sequence} of $\Sigma_{1}^{0}(\mathcal{G})$ sets if there is a total computable function $\phi:\omega^2 \to \omega$ and a computable enumeration $(S_i)$ of $\mathcal{G}$ such that
\begin{equation}
B_n = \cup_m S_{\phi(n,m)}.
\end{equation}
The intersection $\cap_n B_n$ is a $\Pi_{2}^{0}(\mathcal{G})$ set. If the Wiener measure of $B_n$ converges effectively to $0$ as $n\to \infty$, the set $\cap_n B_n$ is said to be a $\Pi_{2}^{0}(\mathcal{G})$ set of constructive measure $0$.

The analogue of Theorem \ref{thm1.1} is the following.
\begin{theorem}\label{thm5.4}\cite{Fouche2}
Let $\mathcal{G}$ be an effectively generated algebra. If $x$ is a complex oscillation, then $x$ is in the complement of every $\Pi_{2}^{0} (\mathcal{G})$ set of constructive measure $0$.
\end{theorem}
We also have an analogue for the effective version of the Borel-Cantelli lemma for Wiener measure.
\begin{theorem}\label{thm5.5}\cite{Fouche2}
If $(A_k)$ is a uniform sequence of $\Sigma_{1}^{0}$ sets with $\sum_k \mathbb{W}(A_k)<\infty$, then for a complex oscillation $x$ there is some $m$ such that $x \notin A_k$ for $k\geq m$.
\end{theorem}

We can now proceed to construct algorithmically random Brownian motion. Suppose we have a sequence of (real) independent random variables $X_0, X_1, Y_1, X_2, Y_2,\dots$ such that each is a normal random variable of mean $0$ and variance $1$. From here on, we will refer to such a sequence simply as a \emph{normal sequence}. We define the Fourier-Wiener series for $t\in [0,1]$ as
\begin{equation}\label{eq5.1}
W(t) = X_0 t +\sqrt{2}\sum_{n=1}^{\infty}\frac{1}{2\pi n}(X_n \sin 2\pi nt +Y_n (1-\cos 2\pi nt)).
\end{equation}
As usual, we suppress the argument $\omega$ in each of the random variables. When the series is parametrised by some $\alpha \in KC$, we denote it by $W_{\alpha}(t)$. By Theorem 2 on p236 of~\cite{Kahane}, the series represents the sample paths of a real Brownian motion on $[0,1]$. In this section we show that the normal sequence can be chosen to be computable and in such a way that the series (\ref{eq5.1}) represents an algorithmically random Brownian motion. The sequence in question was used by Fouch\'{e} in~\cite{Fouche} in order to construct the algoritmically random Brownian motion as a Franklin-Wiener series. As in~\cite{Fouche}, we show that not only does the series converge to a continuous function, but it is a complex oscillation in the sense of Asarin and Pokrovskii~\cite{Asarin}.

We now show how to construct a normal sequence which is parametrised by $\mathcal{N}$, as was first demonstrated in~\cite{Fouche}. Define the function $g:[0,1]\to \mathbb{R}$ by
\begin{equation}
\alpha = \frac{1}{\sqrt{2\pi}}\int_{-\infty}^{g(\alpha)}e^{-t^2 /2} dt, \quad \alpha \in (0,1).
\end{equation} 
Note that $g$ is a computable function. Now, fix a computable bijection $\varphi :\omega^2 \to \omega$. Associate to $\alpha \in \mathcal{N}$ the sequence $(\beta_n)$ such that the $k$th term of $(\beta_n)$ is given by $\beta_k = \alpha_{\varphi (k,0)} \alpha_{\varphi (k,1)}\alpha_{\varphi (k,2)}\dots$. Set $\xi _n (\alpha ) = g(\beta_n)$. The sequence $(\xi_n )$ now forms a normal sequence with respect to Lebesgue measure~\cite{Fouche}. For our later notational convenience, we set
\begin{equation}\label{eq5.11}
X_n = \xi_{2n} \quad \textrm{ and } Y_n = \xi_{2n+1}, \quad n=0,1,2,\dots
\end{equation}

The following theorem, which was already used in proving Theorem \ref{thm4.1}, will be needed to prove Theorem \ref{thm5.2}. The constants have been modified slightly for our purposes.
\begin{theorem}\label{thm5.1}(\cite{Kahane}, p69) Consider the random trigonometric polynomial 
\begin{equation}
P(t) = \sum_{n=1}^{k} \xi_n f_n(t)
\end{equation}
where the $f_n$ are real or complex trigonometric polynomials of degree less than or equal to $N$, $\xi_n$ a normal (or subnormal) sequence and $k<\infty$. Then 
\begin{equation}
\mathbb{P} \left( \| P\|_{\infty} \geq 6(\sum \|f_n \|^{2}_{\infty} \log N)^{\frac{1}{2}}\right) \leq \frac{8\pi}{N^2}.
\end{equation}
\end{theorem}
The particular choice of constants is justified by setting $\rho = 2\pi N^2$ and $\kappa = N^2/4\pi$ in Theorem 1, p68 of~\cite{Kahane}.

\begin{theorem}\label{thm5.2}
Let $\{X_0, X_n,Y_n: n=1,2,\dots \}$ be as defined in (\ref{eq5.11}). Then the series (\ref{eq5.1}) converges to a continuous function whenever $W= W_{\alpha}$ for some $\alpha \in KC$.  
\end{theorem}
\textbf{Proof.} The linear term clearly does not affect convergence, and we shall focus our attention on the sum. Firstly, consider the sum
\begin{equation}\label{eq5.2}
\sum_{n=1}^{\infty} \frac{X_n}{2\pi n}\sin 2\pi nt.
\end{equation}
Note that this is merely a sum of cosines, such as we have been using, with a constant phase shift. The factor $1/2\pi$ will not have an effect on the convergence, and we dispense with it from here on, and also when bounding the terms involving $Y_n$. With $N_k = 2^{2^k}$ as in Theorem \ref{thm4.1}, let
\begin{equation}\label{eq5.5}
P_k (t) = \sum_{n=N_k}^{N_{k+1}-1}\frac{X_n}{n}\sin 2\pi nt,\quad k=0,1,2,\dots
\end{equation}
The series (\ref{eq5.2}) can then be written as $\sum_{k=0}^{\infty} P_k (t)$. Each $P_k$ is a trigonometric polynomial of degree less that $N_{k+1}$. Using Theorem \ref{thm5.1} as in Theorem \ref{thm4.1} (but invoking Theorem \ref{thm5.5} instead of Lemma \ref{lem1.1}), we have that for $\alpha \in KC$ and all large $k$, 
\begin{equation}
\| P_k \|_{\infty} \leq 6\left( 2^{k/2}\left( \sum_{n=2^{2^k}}^{2^{2^{k+1}}-1}\frac{1}{n^2} \right)^{\frac{1}{2}}\right).
\end{equation}
Now, since
\begin{equation}
\sum_{n=2^{2^k}}^{2^{2^{k+1}}-1}\frac{1}{n^2} \leq \left( 2^{2^{k-1}}\right)^{-2},
\end{equation}
we have that
\begin{equation}\label{eq5.3}
\| P_k \|_{\infty} \leq 6 \sqrt{\log 2}\cdot 2^{\frac{k+1}{2}}2^{-2^{k-1}} < 6\cdot 2^{-2^{k-2}} =6N_{k-2}^{-1}.
\end{equation}
The same technique will be used to compute the rate of convergence of the rest of the sum. Set
\begin{equation}\label{eq5.6}
Q_k (t) = \sum_{n=N^k}^{N^{k+1}-1}\frac{Y_n}{n} (1-\cos 2\pi nt).
\end{equation}
Writing $1-\cos 2\pi nt = 2\sin^2 \pi nt$, we observe that Theorem \ref{thm5.1} must now be applied to trigonometric polynomials of degree less than or equal to $2N$. Otherwise, everything else used to obtain the previous estimate (\ref{eq5.3}) still applies, and we obtain
\begin{eqnarray}
\| Q_k \|_{\infty} &\leq & 6\left( \log (2N_{k+1})\sum_{N_{k}}^{N_{k+1}-1} \frac{1}{n^2}\right)^{\frac{1}{2}}\leq 6\sqrt{\log 2} \left( 2^{k+1}+1 \right)^{\frac{1}{2}}\left( \sum_{n=2^{2^k}}^{2^{2^{k+1}}-1}\frac{1}{n^2}\right)^{\frac{1}{2}} \\ & <& 6\cdot 2^{\frac{k+2}{2}}2^{-2^{k-1}} < 6\cdot 2^{-2^{k-2}}= 6N_{k-2}^{-1}.
\end{eqnarray}
for each $\alpha \in KC$. Since $\sum_k N_{k-2}^{-1}$ converges, the series (\ref{eq5.1}) converges uniformly on $[0,1]$. \qed

Using the above estimates, we can now prove
\begin{theorem}
The series (\ref{eq5.1}) converges to a complex oscillation whenever $W=W_{\alpha}$, for each $\alpha \in KC$.
\end{theorem}
The idea of the proof is identical to that of Proposition 2 of~\cite{Fouche}, and can be adapted to the present case using the estimates of Theorem \ref{thm5.1}.

\textbf{Proof.} Suppose that $\alpha \in KC$, but that $W_{\alpha}$ is not a complex oscillation. As in~\cite{Fouche}, we can then write
\begin{equation}
W_{\alpha} \in \bigcap_n \bigcup_m T_{n,m}
\end{equation}
where each $T_{n,m}$ is of the form
\begin{equation}\label{eq5.4}
\bigcap_{i<k} [A_i < X(t_i) <B_i].
\end{equation}
Here, $A_i$, $B_i$ and $t_i$ are dyadic rationals and $k$ an integer, and an effective description of them all can be obtained from $n$ and $m$. Moreover, $\mathbb{W}(\cup_m T_{n,m})$ converges effectively to $0$ as $n\to \infty$. Without loss, we can assume that 
\begin{equation}\label{eq5.9}
\mathbb{W}\left( \bigcup_m T_{n,m}\right)\leq 2^{-n}.
\end{equation} 
Setting $N_j = 2^{2^j}$ as in the previous proof, let
\begin{equation}
W_{\beta}^{N_j}(t) = X_0 t +\sqrt{2}\sum_{n=1}^{N_j -1}\frac{1}{2\pi n}(X_n (\beta) \sin 2\pi nt +Y_n(\beta)(1-\cos 2\pi nt)).
\end{equation}
Let $T_{n,m}$ be as in (\ref{eq5.4}), noting that the parameter $k$ can be effectively retrieved from $n$ and $m$. Define $S_{n,m}$ by
\begin{equation}
\beta \in S_{n,m} \iff \forall_{i<k} \exists_{j> n+m+1}(A_i +N_{j}^{-1/6}< W_{\beta}^{N_j}(t_i) <B_i - N_{j}^{-1/6}).
\end{equation}
The relation $\beta \in S_{n,m}$ is $\Sigma_1^0$ in $\beta$, $n$ and $m$, keeping in mind that $k$ is effectively retrievable from $n$ and $m$. 

By the previous proof, we can bound the tail of the series $P_j +Q_j$, as defined in (\ref{eq5.5}) and (\ref{eq5.6}). When necessary, we shall indicate the dependence on a specific $\beta \in \mathcal{N}$ by writing $P_{j}^{\beta}$ and $Q_{j}^{\beta}$. Specifically, for each $j$ we have that
\begin{equation}\label{eq5.8}
\|P_j(t)+Q_j(t)\|_{\infty} \leq 12 N_{j-2}^{-1}
\end{equation} for $k$ large enough. 
Therefore, we get the bound
\begin{equation}
\frac{1}{2\pi}\left\| \sum_{i=j}^{\infty}(P_i(t)+Q_i(t))\right\|_{\infty} \leq \frac{12}{\pi} 2^{-2^{j-2}}= \frac{12}{\pi} N_{j-2}^{-1} 
\end{equation} for large $k$.
This allows us to state that, for a given $\alpha \in KC$, we can find some $j_0 = j_0 (\alpha)$ such that 
\begin{equation}\label{eq5.7}
|W_{\alpha}(t) - W_{\alpha}^{N_{j}}(t)|\leq \frac{1}{2\pi}\left\| \sum_{i=j}^{\infty}(P_i(t)+Q_i(t))\right\|_{\infty} \leq \frac{12}{\pi} N_{j-2}^{-1}.
\end{equation}
for all $j\geq j_0$. For $\alpha \in KC$ specified and $j_0$ determined, we want to show that 
\begin{equation}
\alpha \in \bigcap_{n\geq j_0}\bigcup_{m} S_{n,m}.
\end{equation}

For convenience, we now set $C= 12/\pi$. Let $n\geq j_0$, and choose $m$ such that $x_{\alpha}\in T_{n,m}$. We can now find some $L>0$ such that for all $i<k$,
\begin{equation}
A_i +\frac{1}{L}< W_{\alpha}(t_i) <B_i -\frac{1}{L}.
\end{equation}
Now choose $j> n+m+1$ such that $CN_{j-2}^{-1}+N_{j}^{-1/6}<L^{-1}$. By (\ref{eq5.7}) and the fact that $j>j_0$, 
\begin{equation}
W_{\alpha}(t_i) -CN_{j-2}^{-1} < W_{\alpha}^{N_j}(t_i) < W_{\alpha}(t_i) +CN_{j-2}^{-1}.
\end{equation}
From the choice of $j$, 
\begin{equation}
A_i +N_{j}^{-1/6} < W_{\alpha}^{N_j}(t_i) <B_i - N_{j}^{-1/6},
\end{equation}
and therefore $\alpha \in S_{n,m}$.

The aim is to show that $\cup_m S_{n,m}$ tends to $0$ effectively, which would mean that $\alpha$ is contained in a $\Pi_2^0$ set of constructive measure $0$, in contradiction to the assumption that $\alpha \in KC$. Therefore, suppose that $\beta \in S_{n,m}$, where $n\geq j_0$. We will now establish that if $\| P_{j}^{\beta} +Q_{j}^{\beta} \|_{\infty} \leq 12N_{j-2}^{-1}$ for all $j\geq n+m+1$, then $W_\beta \in T_{n,m}$, where $P_{j}^{\beta}$ and $Q_{j}^{\beta}$ are the polynomials (\ref{eq5.5}) and (\ref{eq5.6}) for the variables $X_n$ and $Y_n$ parametrised by $\beta$.
By definition of $S_{n,m}$, there are some $i$ and $j>n+m+1$ such that
\begin{equation}
A_i + N_{j}^{-1/6} < W_{\beta}^{N_j}(t_i) <B_i - N_{j}^{-1/6}. 
\end{equation}
Because of the assumed bounds on $P_j +Q_j$, for large $j$, 
\begin{equation}
|W_{\beta}(t_i) - W_{\beta}^{N_{j}}(t_i)|\leq CN_{j-2}^{-1},
\end{equation}
and so
\begin{equation}
A_i + N_{j}^{-1/6} - CN_{j-2}^{-1} \leq W_{\beta}(t_i) \leq B_i - N_{j}^{-1/6} + CN_{j-2}^{-1}.
\end{equation}
Except for a few initial terms, $N_{j}^{-1/6} - CN_{j-2}^{-1}>0$, and we can assume that $j_0$ is also large enough for this to be true. For $n\geq j_0$ then, $W_{\beta} \in T_{n,m}$. This will also hold under the assumption that $\beta \in S_{n,m}$ and that both $\|P^{\beta}_{j}\|_{\infty}$ and $\|Q^{\beta}_{j}\|_{\infty}$ are bounded from above by $6N_{j-2}^{-1}$ for large $j$, since then the condition on $\| P_{j}^{\beta} +Q_{j}^{\beta} \|_{\infty}$ is clearly implied. If we now define $U_{n,m}$ by the relation
\begin{equation}
\beta \in U_{n,m} \iff \exists_{j>n+m+1}\left( (\|P_{j}^{\beta}\|_{\infty} > 6N_{j-2}^{-1}) \lor (\| Q_{j}^{\beta}\|_{\infty}> 6N_{j-2}^{-1})\right),
\end{equation} 
we can state that
\begin{equation}
\left( \beta \in S_{n,m} \right) \implies \exists_m [(W_{\beta}\in T_{n,m})\lor \beta \in U_{n,m}].
\end{equation}
Utilising Theorem \ref{thm5.1} again, we see that 
\begin{equation}
\mathbb{P} [\|P_{j}^{\beta}\|_{\infty}>6N_{j-2}^{-1}] < \mathbb{P}\left[ \|P_{j}^{\beta}\|_{\infty} > 6\left( 2^{j+1}\log 2 \sum_{N_j}^{N_{j+1}-1} n^{-2}\right)^{1/2} \right] <8\pi/N_{j+1}^{2},
\end{equation} and the same for $Q_{j}^{\beta}$.

The final component necessary to bound the measure of $S_{n,m}$ is to note that, for $n\geq j_0$,
\begin{equation}
\mathbb{P} \left( \cup_m [\alpha : W_{\alpha}\in T_{n,m}]\right) = \mathbb{P} \left( [\alpha: W_{\alpha}\in \cup_m T_{n,m}]\right) = 
\mathbb{W}\left( \cup_m T_{n,m} \right)\leq 2^{-n},
\end{equation}
by (\ref{eq5.9}). We now have that 
\begin{equation}
\mathbb{P} \left( \cup_m S_{n,m} \right) \leq  2^{-n}+\sum_{j>m+n+1}\frac{16\pi}{N_{j+1}^{2}},
\end{equation}
which converges to $0$ effectively as $n\to \infty$. This contradicts the assumption that $\alpha \in KC$, which implies that $W_{\alpha}$ is indeed a complex oscillation. \qed


\begin{thebibliography}{9}

\bibitem{Asarin}Asarin, E.A. and Pokrovskii, A.V., Use of the Kolmogorov complexity in analyzing control system dynamics. 
 \emph{Automation and Remote Control}, 47:21--28, 1986.

\bibitem{Davie} Davie, G., The Borel-Cantelli lemmas, probability laws and Kolmogorov complexity. \textit{Ann. Prob.} 29 (4), 1426--1434.

\bibitem{Fouche2} Fouch{\'e}, W.L., Arithmetical representations of Brownian motion. \textit{J. Symbolic Logic} 65 421--442, 2000.

\bibitem{Fouche} Fouch{\'e}, W.L., The descriptive complexity of Brownian motion. \textit{Advances in Mathematics}. 155(2):317--343, 2000.

\bibitem{Fouche3} Fouch{\'e}, W.L., Descriptive complexity and reflective properties of combinatorial configurations. \textit{J. London Math. Soc.} 54, 199--208, 1996. 

\bibitem{FMD} Fouch{\'e}, W.L., Mukeru, S., Davie, G., Fourier spectra of measures associated with algorithmically random Brownian motion. \textit{Log. Methods Comput. Sci.}, Vol. 10 (3:20) 1--24, 2014. 


\bibitem{Hanssen} Kjos-Hanssen, B., Nerode, A., The law of the iterated logarithm for algorithmically random Brownian motion, in: Proceedings on Logical Foundations of Computer Science, LFCS 2007, in: Lecture Notes in Computer Science, vol. 4514, 2007, pp 310--317.


\bibitem{Kahane} Kahane, J-P., \emph{Some Random Series of Functions}, 2nd ed.
Cambridge University Press, 1985

\bibitem{MartinLof} Martin-L\"{o}f, P., The definition of random sequences. \textit{Inform. and Control}, 9 602--619, 1966.

\bibitem{Zygmund} Zygmund, A., \emph{Trigonometric Series}, Cambridge University Press, 1959, vol. 1.

\end{thebibliography}
\end{document}